\documentclass[centertags,12pt]{amsart} \usepackage{latexsym} 
\usepackage{amssymb} \usepackage{graphics}

\numberwithin{equation}{section}

\makeatletter \renewcommand{\subsection}{\@startsection 
{subsection}{2}{0mm}{\baselineskip}{-0.25cm} 
{\normalfont\normalsize\em}} \makeatother

\newtheorem{theorem}{Theorem}[section] 
\newtheorem{proposition}[theorem]{Proposition} 
\newtheorem{corollary}[theorem]{Corollary} 
\newtheorem{lemma}[theorem]{Lemma}

{\theoremstyle{definition} 
 
\newtheorem{example}[theorem]{Example} 
 
}

\theoremstyle{remark} \newtheorem{remark}[theorem]{Remark}

\def\bF{\mathbb F} \def\bN{\mathbb N} \def\bP{\mathbb P} 
\def\bZ{\mathbb Z}

\def\1{\mathbf 1} \def\a{\mathbf a} \def\b{\mathbf b} \def\c{\mathbf 
c}    
  \def\cero{\mathbf 0}

   \def\K{\mathbf 
K}  \def\R{\mathbf R}

\def\cB{\mathcal B}   
   
\def\cL{\mathcal L} \def\cM{\mathcal M} \def\cO{\mathcal O} 
 \def\cX{\mathcal X}  
\def\cU{\mathcal U}

\def\gaps{{\rm Gaps}} \def\supp{{\rm Supp}}  
 \def\dim{{\rm dim}}

\begin{document}

\author[C. Munuera]{Carlos Munuera} \author[F. Torres]{Fernando 
Torres}

\title[Algebras admitting well agreeing n-weights]{The structure of 
algebras admitting \\ well agreeing near weights}

\address{Department of Applied Mathematics, University of Valladolid 
(ETS Arquitectura) 47014 Valladolid, Castilla, Spain.}

\email{cmunuera@modulor.arq.uva.es}

\address{IMECC-UNICAMP, Cx.P. 6065, 13083-970, Campinas SP-Brazil.}

\email{ftorres@ime.unicamp.br}


\thanks{{\em 2000 Math. Subj. Class.}: Primary 94, Secondary 14}

\thanks{The authors were supported respectively by the ``Junta de 
Castilla y Le\'on", Espa\~na, under Grant VA020-02, and CNPq-Brazil 
(306676/03-6) and PRONEX (66.2408/96-9)}


   \begin{abstract} We characterize algebras admitting two well 
agreeing near weights $\rho$ and $\sigma$. We show that such an 
algebra $R$ is an integral domain whose quotient field $\K$ is an 
algebraic function field of one variable. It contains two 
places $P, Q\in {\mathbb P}(\K)$ such that $\rho$ and $\sigma$ are 
derived from the valuations associated to $P$ and $Q$. Furthermore 
$\bar R=\cap_{S\in \bP(\K)\setminus\{P,Q\}} \cO_S.$
   \end{abstract}

\maketitle

   \section{Introduction}\label{s1}
   
Algebraic Geometric codes (or AG codes, for short) were constructed by 
Goppa \cite{goppa1}, \cite{goppa2}, based on a curve ${\mathcal X}$ 
over a finite field ${\mathbb F}$ and two rational divisors $D$ and 
$G$ on ${\mathcal X}$, where $D$ is a sum of pairwise distinct points 
and $G=\alpha_1P_1+\dots+\alpha_mP_m$, with $P_i\not\in\supp D$. Soon 
after its introduction, AG codes became a very important tool in 
Coding Theory; for example, Tsfasman, Vladut and Zink \cite{tsf-vla-z} 
showed that the Varshamov-Gilbert bound can be attained by using these 
codes.  However, the study of AG codes relies on the use of algebraic 
geometric tools, which is difficult for non specialists in Algebraic 
Geometry. In 1998, H\o holdt, Pellikaan and van Lint presented a 
construction of AG codes `without Algebraic Geometry'; that is, by 
using elementary methods only \cite{h-vl-p} (see also \cite{feng}). 
These methods include order and weight functions over an ${\mathbb 
F}$-algebra and Semigroup Theory mainly. From that paper, order 
domains and order functions and the corresponding obtained codes have 
been studied by many authors; to mention a few of them: Pellikaan 
\cite{ord2}, Geil and Pellikaan \cite{geil-p} and Matsumoto 
\cite{matsumoto}.

The approach given by H\o holdt, Pellikaan and van Lint allows us to 
do with the so called `one point' AG codes; that is, when the divisor 
$G$ is a multiple of a single point, $G=\alpha P$. A generalization of 
the same idea to arbitrary AG codes ($m\geq 1$) was given in 
\cite{norders}. To that end, variations of order and weight functions 
over an ${\mathbb F}$-algebra $R$ --the so called {\em near order} and 
{\em near weight functions}-- are introduced.

In the present paper, we characterize algebras $R$ admitting two well 
agreeing near weights (see Section \ref{s2} for explanation of this 
concept), $\rho$ and $\sigma$, as being certain subalgebras of the 
regular function ring of an affine variety of type $\cX\setminus 
\{P,Q\}$, where $\cX$ is a projective, geometrically irreducible, 
non-singular algebraic curve and $P$ and $Q$ are two different points 
of $\cX$. We will also show that $\rho$ and $\sigma$ are defined by 
the valuations at $P$ and $Q$ respectively (see Theorem \ref{thm5.1} 
in Section \ref{s5}). This result is essentially analogous to the 
characterization of algebras admitting a weight function given by 
Matsumoto \cite{matsumoto}.

For simplicity, throughout this paper we shall use the language in 
terms of algebraic function fields instead of algebraic curves.

   \section{Near weights}\label{s2}
   
In this section we recall the concept of near weight and discuss some of 
its properties. Throughout, let $R$ be an algebra over a field $\bF$. 
We always assume that $R$ is commutative and $\bF \subsetneq R$. For a 
function $\rho:R\to \bN_0\cup\{-\infty\}$, let us consider the sets
  \begin{align*} \cU^* =\cU^*_\rho & := \{r\in R\setminus\{0\}: 
\rho(r)\leq \rho (1)\};\\ \cM =\cM_\rho & := \{r\in R: 
\rho(r)>\rho(1)\}
  \end{align*} 
and $\cU=\cU_\rho:=\cU^*\cup\{0\}$. The function $\rho$ 
is called a {\em near weight} (or a {\em n-weight}, for short) if the 
following conditions are satisfied. Let $f, g, h\in R$;
  \begin{enumerate} \item[(N0)] $\rho(f)=-\infty$ if and only if 
$f=0$; \item[(N1)] $\rho(\lambda f)=\rho(f)$ for $\lambda\in\bF^*:= 
\bF\setminus\{0\}$; \item[(N2)] $\rho(f+g)\leq \rho(f)+\rho(g)$; 
\item[(N3)] If $\rho(f)<\rho(g)$, then $\rho(fh)\leq \rho(gh)$. 
Furthermore, if $h\in\cM$ then $\rho(fh)<\rho(gh)$; \item[(N4)] If 
$\rho(f)=\rho(g)$ with $f,g\in\cM$, then there exists 
$\lambda\in\bF^*$ such that $\rho(f-\lambda g)<\rho(f)$; \item[(N5)] 
$\rho(fg)\leq \rho(f)+\rho(g)$ and equality holds if $f,g\in\cM$.
  \end{enumerate} Near weights were introduced in \cite{norders} in 
connection with an elementary construction of algebraic geometric 
codes.  After a normalization, we can assume $\rho(f)=0$ for $f\in 
\cU^*$ and $\gcd \{ \rho(f) : f\in \cM\}=1$; see \cite[Sect. 3.2]{norders}. A n-weight becomes a {\em 
weight function}, as defined in H\o holdt, van Lint and Pellikaan if 
and only if $\cU=\bF$ \cite[Lemma 3.3]{norders}.

For given two n-weights $\rho$ and $\sigma$ over the $\bF$-algebra 
$R$, set
  $$ H=H(R):=\{(\rho(f),\sigma(f): f\in R^*\}\, ,
  $$ where $R^*=R\setminus\{0\}$. We say that $\rho$ and $\sigma$ {\em 
agree well} if $\# (\bN^2\setminus H)$ is finite and 
$\cU_\rho\cap\cU_\sigma=\bF$. In the next section we will prove that 
$H$ is a semigroup so that this definition will in fact be compatible 
with the one given in \cite{norders}. As said before, our purpose in 
this paper is to characterize the algebras $R$ admitting well agreeing 
n-weights. These algebras exist, as the next example shows.
  \begin{example}\label{ex2.1} Let ${\mathbf K}$ be an algebraic 
function field of one variable over $\mathbb F$, such that $\mathbb F$ 
is the full constant field of $\mathbf K$. For a place $S$ of $\mathbf 
K$, let $\cO_S$ be the local ring at $S$ and $v_S$ its corresponding 
valuation. Let $P, Q$ be two different places of $\mathbf K$. 
Consider an $\bF$-algebra $R\subseteq\mathbf K$ and 
define
$$ \varrho(f):=\begin{cases} -\infty & \text{if $f=0$},\\ 0 & 
\text{if $v_P(f)\geq 0$},\\ -v_P(f) & \text{if $v_P(f)<0$},
 \end{cases} \quad\text{and}\quad \varsigma(f):=\begin{cases} -\infty 
& \text{if $f=0$},\\ 0 & \text{if $v_Q(f)\geq 0$},\\ -v_Q(f) & 
\text{if $v_Q(f)<0$}\, .
 \end{cases}
 $$ 
If $R=R(P,Q):=\bigcap_{S\neq P,Q}\cO_S$, then $\varrho$ and 
$\varsigma$ are well agreeing n-weights over 
$R$.
  \end{example} To end this section, we state a property of n-weights 
that we shall need later.
  \begin{lemma}\label{lemma4.3} Let $f,g\in R^*$ such that $\rho(f)>0$ 
and $\rho(g)=0.$ Then there exists $\lambda\in\bF$ such that 
$\rho(f(g-\lambda))<\rho(f).$
  \end{lemma}
  \begin{proof} According to (N5), $\rho(fg)\leq \rho(f)$. If 
$\rho(fg)= \rho(f)$, by (N4) there exists $\lambda\in\bF$ such that 
$\rho(f(g-\lambda))=\rho(fg-\lambda f)<\rho(f)$.
  \end{proof}

   \section{The semigroup structure}\label{s3}

Let $\rho$ and $\sigma$ be two well agreeing n-weights defined on an 
$\bF$-algebra $R$. We generalize the definition of the set $H=H(R)$ 
stated in Section \ref{s2} to any $S\subseteq R$ by setting
   $$ H(S):=\{(\rho(f),\sigma(f)): f\in S^*\}\subseteq \bN_0^2\, ,
   $$ where $S^*:=S\setminus \{ 0\}$. We shall see that $H$ is a 
semigroup. To that end, we need some preliminary results. For given 
$\a=(a_1,a_2)$ and $\b=(b_1,b_2)$ elements of $\bN_0^2$, the {\em 
least upper bound} of $\a$ and $\b$ is defined as (cf. 
\cite{gretchen1}, \cite{gretchen2})
   $$ {\rm lub}(\a,\b):=(\max\{a_1,b_1\},\max\{a_2,b_2\})\, .
   $$
   \begin{lemma}\label{lemma3.1} Let $f,g\in R^*$. Set 
$\a:=(\rho(f),\sigma(f))$ and $\b:=(\rho(g),\sigma(g)).$ Then there 
exist $\lambda,\mu\in\{0,1\}$ such that
  $$ {\rm lub}(\a,\b)=(\rho(\lambda f+\mu g),\sigma(\lambda f+\mu 
g))\, .
  $$ In particular, if $f,g \in S\subseteq R$ and $S$ is closed 
under sum, then ${\rm lub}(\a,\b)\in H(S).$
   \end{lemma}
   \begin{proof} If $\a=\b$ the result is obvious. Otherwise, we can 
assume $\rho(f)<\rho(g)$. If $\sigma(f)\leq \sigma(g)$, then ${\rm 
lub}(\a,\b)=\b$. On the contrary, if $\sigma(f)> \sigma(g)$ then
  $$ \max\{\rho(f),\rho(g)\} = \rho(f+g)\, \quad\text{and}\quad 
\max\{\sigma(f),\sigma(g)\}= \sigma(f+g)
  $$ and hence ${\rm lub}(\a,\b)=(\rho(f+g),\sigma(f+g))$. The second 
part of the lemma is clear.
    \end{proof}
   \begin{proposition}\label{prop3.1} Let $S\subseteq R$ be a closed 
subset under sum and product. Then $H(S)$ is closed for the sum$;$ 
that is$,$ if $\a,\b\in H(S),$ then $\a+\b\in H(S).$
   \end{proposition}
   \begin{proof} Let $\a=(\rho(f),\sigma(f))$ and 
$\b=(\rho(g),\sigma(g))$ with $f,g\in S^*$. If $\a=\cero$, the result is 
clear. If the integers $\rho(f),\sigma(f),\rho(g),\sigma(g)$ are all 
positive; that is, $f,g\in\cM_\rho\cap\cM_\sigma$, the result follows 
from property (N5) of n-weights. Then assume $\rho(f)>0$ and 
$\sigma(f)=0$. There are three possibilities:
  \begin{itemize} \item If $\rho(g)=0$ and $\sigma(g)>0$, then 
$\a+\b={\rm lub}(\a,\b)\in H(S)$ according to Lemma \ref{lemma3.1}; 
\item If $\rho(g)>0$ and $\sigma(g)=0$, then 
$\a+\b=(\rho(fg),\sigma(fg)) \in H(S)$ by (N5); \item If $\rho(g)>0$ 
and $\sigma(g)>0$, then $\rho(fg)=\rho(f)+\rho(g)$ and $\sigma(fg)\leq 
\sigma(g)$ and hence $\a+\b={\rm lub}(\a,\c)\in H(S)$, where 
$\c=(\rho(fg),\sigma(fg))$.
  \end{itemize}
  \end{proof}
  \begin{corollary}\label{cor3.1} Let $R'$ be a $\bF$-subalgebra of 
$R.$ Then $H(R')$ is a semigroup$.$
  \end{corollary} Next we consider the following sets associated to 
the semigroup $H=H(R)$:
   $$ H_x:=\{(m,0)\in H \}\, ,\quad H_y:=\{(0,n)\in H\}\, ,
   $$ and their projections
   $$ \bar H_x:=\{ m : (m,0)\in H\}\, , \quad \bar H_y:=\{ n : 
(0,n)\in H\}\, .
   $$ Clearly $\bar H_x$ and $\bar H_y$ are numerical semigroups of 
finite genus. For $n\in\bN_0$, set
  $$ x_H(n) := \min\{m\in\bN_0: (m,n)\in H\}\quad\text{and}\quad 
y_H(n) := \min\{m\in\bN_0: (n,m)\in H\}\, .
  $$
   \begin{lemma}\label{lemma3.2} If $y_H(n)>0,$ then 
$x_H(y_H(n))=n>0.$
   \end{lemma}
   \begin{proof} Let $f\in R^*$ such that $\rho(f)=n$ and 
$\sigma(f)=y_H(n)$. By definition, $x_H(y_H(n))\leq n$. If $\rho(g)<n$ 
and $\sigma(g)=y_H(n)$ for some $g\in R$, then there exists 
$\lambda\in\bF$ such that $\sigma(f-\lambda g)<y_H(n)$. Since 
$\rho(f-\lambda g)=\rho(f)$, this is a contradiction.
    \end{proof}
    \begin{corollary}\label{cor3.2} {\rm (cf. \cite{kim}, \cite[Cor. 
4.8]{norders})} It holds that $n\in\gaps(\bar H_x)$ if and only if 
$y_H(n)\in\gaps(\bar H_y).$ In particular$,$ the semigroups $\bar H_x$ 
and $\bar H_y$ have equal genus$.$
    \end{corollary} We consider now the following subsets of $H$:
  \begin{align*} \tilde\Gamma=\tilde\Gamma(H) & := \{(m,y_H(m)): 
m\in\gaps(\bar H_x)\}=\{(x_H(n),n): n\in\gaps(\bar H_y)\}, \\ 
\Gamma=\Gamma(H) & :=\{(m,y_H(m)), (x_H(m),m) : m\in\bN_0\}=\tilde 
\Gamma\cup H_x \cup H_y.
   \end{align*} Note that $\tilde\Gamma$ is well defined according to 
Lemma \ref{lemma3.2}. The result below allows a nice description of 
the semigroup $H$.
   \begin{proposition}\label{prop3.2} {\rm (cf. \cite{kim}, 
\cite{gretchen1})}
  $$ H=\{{\rm lub}(\a,\b): \a, \b\in \Gamma\}\, .
  $$
   \end{proposition}
   \begin{proof} According to Lemma \ref{lemma3.1}, ${\rm lub}(\a, 
\b)\in H$ for all $\a,\b\in H$. Conversely, each $\a=(a_1,a_2)$ can be 
written as $\a={\rm lub}((a_1,y_H(a_1)), (x_H(a_2),a_2))$.
   \end{proof} For every $\a\in H$ take an element $\phi_\a\in R^*$ 
such that $(\rho(\phi_\a),\sigma(\phi_\a))=\a$, and set
  $$ \cB:=\{\phi_\a: \a\in\Gamma\}\, .
  $$
   \begin{proposition}\label{prop3.3} The set $\cB$ is a basis of $R$ 
as a $\bF$-vector space$.$
   \end{proposition}
   \begin{proof} Since every two points $\a\neq\b\in\Gamma$ lie in 
different row and column, the set $\cB$ is linearly independent, 
according to property (N2) of n-weights. To see that $\cB$ generate 
$R$ take an element $f\in R^*$. Let us assume first that $\sigma(f)=0$ 
and use induction on $\rho(f)$. If $\rho(f)=0$ the result follows from 
the fact that $\cU_\rho\cap\cU_\sigma=\bF$. If $\rho(f)=k>0$, take 
$\phi_\a\in\Gamma$ with $\a=(k,0)$. There exists $\lambda\in\bF$ such 
that either $\lambda\phi_\a=f$ or $\rho(f-\lambda\phi_\a)<k$ and 
$\sigma(f-\lambda\phi_\a)=0$. By induction hypothesis, all elements 
$g$ with $\sigma(g)=0$ and $\rho(g)<k$ are generated by $\cB$ and 
hence $f$ is generated by $\cB$. According to Lemma \ref{lemma3.1} and 
Proposition \ref{prop3.2}, the general case $\sigma(f)>0$ follows now 
by induction on $\sigma(f)$.
   \end{proof} For $(m,n)\in\bN_0^2$ write
  $$ \Delta(m,n):=\{(m,\ell): \ell<n\}\cup\{(\ell,n): \ell<m \}\, .
   $$ Let $\gaps(H)$ denotes the set of gaps of $H$.
   \begin{corollary}\label{cor3.3} {\rm (cf. \cite{ct})} We have
  $$ \gaps(H)=\bigcup_{\a\in\tilde\Gamma} \Delta(\a)\, .
  $$
   \end{corollary}
   \begin{proof} If $(m,n)\in\Delta(\a)$ for some 
$\a\in\tilde{\Gamma}$, then $m<x_H(n)=a_1$ or $n<y_H(m)=a_2$; hence 
$(m,n)\not\in H$. If $(m,n)\not\in\Delta(\a)$ for every 
$\a\in\tilde{\Gamma}$, then $n\geq y_H(m)$ and $m\geq x_H(n)$ and 
hence $(m,n)={\rm lub}((m,y_H(m)),(x_H(n),n))\in H$.
   \end{proof}
   \begin{remark}\label{rem3.1} In the case of Example \ref{ex2.1}, 
$H$ is the Weierstrass semigroup at $P$ and $Q$. A point $(m,n)\in 
{\mathbb N}_0^2$, is a gap of $H$ if and only if 
$\ell(mP+nQ)=\ell((m-1)P+nQ)$ or $\ell(mP+nQ)=\ell(mP+(n-1)Q)$. Homma 
and Kim \cite{homma-kim} noticed that AG codes associated to gaps 
$(m,n)$ where both equalities above hold true, have quite good 
parameters; such gaps are called {\em pure}. Let $\gaps_0(H)$ denotes 
the set of pure gaps of $H$. Then
   $$ 
\gaps_0(H)=\bigcup_{\a\neq\b\in\tilde\Gamma}(\Delta(\a)\cap\Delta(\b))\, 
.
  $$
  \end{remark}
   \begin{remark}\label{rem3.2} For $m,n\in\bN_0$, we can consider the 
subset of $R$
  $$ R(m,n):=\{f\in R:\text{$\rho(f)\leq m$ and $\sigma(f)\leq n$}\}\, 
.
  $$ In \cite{norders}, subsets of this form were used to construct 
codes. Clearly $H(R(m,n))=H(m,n)=\{\a=(a_1,a_2)\in H: \text{$a_1\leq 
m$ and $a_2\leq n$}\}$. Again in the case of Example \ref{ex2.1}, if 
$R=R(P,Q)$, then $H(m,n)=\cL(mP+nQ)$. As a consequence of the 
Proposition \ref{prop3.3}, the set $\{\phi_\a : \a\in\Gamma\cap 
H(m,n)\}$ is a basis of $H(m,n)$.
   \end{remark}

   \section{The structure of the algebra $R$}\label{s4}
   
By keeping the notation of the previous sections, let $R$ be an 
$\bF$-algebra and $\rho$ and $\sigma$ two well agreeing n-weights on 
$R$. The semigroups $\bar H_x$ and $\bar H_y$ are finitely generated 
since they have finite genus. Write
  $$ \bar H_x=\langle m_1,\dots,m_r\rangle\, ,\quad\text{and}\quad 
\bar H_y=\langle n_1,\dots,n_s \rangle
  $$ and define
   $$ \Gamma^+=\Gamma^+(H):=\tilde\Gamma\cup \{ 
(m_1,0),\dots,(m_r,0),(0,n_1),\dots,(0,n_s)\}\subseteq H\, .
   $$
   \begin{lemma}\label{lemma4.1} Let $R'=\bF[\{\phi_\a : 
\a\in\Gamma^+\}]\subseteq R.$ Then $H(R')=H(R).$
   \end{lemma}
   \begin{proof} Clearly $H(R')\subseteq H(R)$. To see the equality, 
according to Proposition \ref{prop3.2} and Lemma \ref{lemma3.1}, it 
suffices to show that $\Gamma\subseteq H(R')$. Let $(m,0)\in H_x$. 
There exist $\alpha_1,\dots,\alpha_r\in\bN_0$ such that $\alpha=\sum 
\alpha_i m_i$. Thus the element
  $$ \phi=\prod \phi_{(m_i,0)}^{\alpha_i}
   $$ belongs to $R'$ and verifies
   $$ \rho(\phi) = \sum \alpha_i\rho(\phi_{(m_i,0)})=\sum\alpha_i 
m_i=m\quad\text{and}\quad \sigma(\phi)\leq \max \{ 
\sigma(\phi_{(m_i,0)}) \}=0
   $$ (because $m_i>0$ and hence $\phi_{(m_i,0)}\in\cM_{\rho}$). Then 
$(m,0)\in H(\R')$. Analogously $H_y\subseteq H(R')$.
    \end{proof}
     \begin{lemma}\label{lemma4.2} Let $R'$ be a $\bF$-subalgebra of 
$R$. If $H(R')=H(R),$ then $R'=R.$
    \end{lemma}
    \begin{proof} Take $f\in\ R$ and let us see that $f\in\ R'$. We 
first consider the case $\sigma(f)=0$. Let us write $\bar H_x=\{ 
\ell_0=0<\ell_1<\ell_2<\dots \}$ and proceed by induction on 
$\rho(f)$. If $\rho(f)=0$ then $f\in\cU_{\rho}\cap\cU_{\sigma}=\bF$ 
and hence $f\in\bF\subseteq R'$. By induction hypothesis assume that 
$f\in R'$ whenever $\rho(f)< \ell_{k+1}$, $k>0$. If $\rho(f)=\ell_k$, 
take $f'\in R'$ such that $\rho(f')=\ell_k$ and $\sigma(f')=0$. Thus, 
there exists $\lambda\in\bF$ such that $\rho(f-\lambda f')<\ell_k$. 
Since $\sigma(f-\lambda f')=0$, we get $f-\lambda f'\in\ R'$ and thus 
$f\in R'$.

Let us prove now the general case by induction on $\sigma(f)$. Assume 
the result true when $\sigma(f)<k+1$. If $\sigma (f)=k$ take $f''\in\ 
R'$ such that $\sigma(f'')=k$. Again there exists $\lambda\in\bF$ such 
that $\sigma(f-\lambda f'')<k$; hence, by induction hypothesis, 
$f-\lambda f''\in R'$ and so $f\in R'$.
    \end{proof}
    \begin{theorem}\label{thm4.1} The $\bF$-algebra $R$ is finitely 
generated over $\bF,$ namely
  $$ R=\bF[\{\phi_\a : \a\in\Gamma^+\}]\, .
  $$
   \end{theorem}
   \begin{proof} It is a direct consequence of Lemmas \ref{lemma4.1} 
and \ref{lemma4.2}.
    \end{proof}
    \begin{proposition}\label{prop4.1} The $\bF$-algebra $R$ is an 
integral domain$.$
   \end{proposition}
   \begin{proof} By \cite[Lemma 3.4]{norders} the set of zero divisors 
of $R$ is contained in $\cU_\rho\cap\cU_\sigma=\bF$; the proof now 
follows as $\rho$ and $\sigma$ are well agreeing by hypothesis.
   \end{proof} In particular, $R$ is isomorphic to an affine 
$\bF$-algebra,
   $$ R\cong \bF[X_1,\dots,X_n]/{\mathfrak q}\, ,
   $$ where ${\mathfrak q}$ is a prime ideal. As an integral domain, 
$R$ admits a field of quotients which we denote by $\K$.
  \begin{theorem}\label{thm4.2} The transcendence degree of $\K$ over 
$\bF$ is one$.$
  \end{theorem} In order to prove this theorem, we need some auxiliary 
results.
  \begin{lemma}\label{lemma4.4} Let $f\in R^*$ and $I=(f)$ be the 
ideal generated by $f.$ The sets $H_x\cap (\bN_0^2\setminus H(I))$ and 
$H_y\cap (\bN_0^2\setminus H(I))$ are both finite$.$
  \end{lemma}
  \begin{proof} If $f\in\bF$ there is nothing to prove. In other case, 
by applying iteratively Lemma \ref{lemma4.3}, there exists $g\in\R^*$ 
such that $\rho(fg)=0$ and hence $\sigma(fg)>0$. Let $\ell_{\sigma}$ 
be the largest gap of $\bar H_y$. Then, for all 
$m>\sigma(fg)+\ell_\sigma$ it holds that $\a=(0,m)\in H(I)$. Indeed, 
let a $\phi\in\R$ be a function such that 
$(\rho(\phi),\sigma(\phi))=(0,m-\sigma(fg))$; then $fg\phi\in I$ and 
$(\rho(fg \phi),\sigma(fg\phi))=\a$. The proof for $H_x$ is analogous.
  \end{proof}
  \begin{proposition}\label{prop4.2} Let $I\subseteq R$ be a proper 
ideal of $R$. Then$,$ as a vector space over $\bF,$ $\dim_\bF(R/I)\leq 
\# \{ \a\in\Gamma : \a\notin H(I)\}.$ In particular$,$ this dimension 
is finite$.$
  \end{proposition}
  \begin{proof} Let $f\in I$, $f\neq 0$, and let $J=(f)$. For every 
$\a\in\Gamma$ take an element $\phi_\a\in \cB$; that is, $(\rho( 
\phi_\a),\sigma(\phi_\a))=\a$. If $\a\in H(J)$ (resp. $\a\in H(I)$) 
take $\phi_\a\in J$ (resp. $\phi_\a\in I$). As we have seen in 
Proposition \ref{prop3.3}, the set $\cB$ is a basis of $R$; hence the 
set of residual classes $\{\phi_\a+I:\a\in\Gamma\}$ form a generator 
system of $R/I$. Now, according to Lemma \ref{lemma4.4} only finitely 
many of these classes are not in $J\subseteq I$.
  \end{proof} {\bf Proof of Theorem \ref{thm4.2}.} According to 
Theorem \ref{thm4.1}, the $\bF$-algebra $R$ is finitely generated over 
$\bF$. Thus the transcendence degree of $\K$ over $\bF$ is equal to the 
Krull dimension of $R$; see Eisenbud \cite[Thm. A p.221]{eisenbud} or 
Matsumura \cite[Ch. 5, Sect. 14]{matsumura}. Take $f\in R^*$ such that 
$f$ is 
not invertible. Such an $f$ exists: it is enough to take $f\in R\setminus 
\bF$. Let $\mathfrak p$ be a minimal prime ideal containing $f$. Then 
${\rm height}(\mathfrak p)=1$ by Krull's Hauptidealsatz; see 
\cite[Thm. 10.2]{eisenbud}. Since (see e.g. 
\cite[Cor. 13.4]{eisenbud} or \cite[Thm. 14.H]{matsumura}),
  $$ 
{\rm height}(\mathfrak p)+\dim(R/\mathfrak p)=\dim(R)\, ,
  $$ 
where `$\dim$' means Krull dimension, and $\dim(R/\mathfrak p)=0$ 
according to Proposition \ref{prop4.2}, we get $\dim(R)=1$.
   \begin{remark}\label{rem4.1} Let $f\in R^*$ and $I=(f)$ be the 
ideal generated by $f$. Let $\a$ and $\b$ be two different points in 
$\Gamma$ and $\phi_\a,\phi_\b\in\cB$. Note that $\phi_\a-\phi_\b\in I$ 
implies ${\rm lub} (\a,\b)\in H(I)$ (as the points $\a,\b$ lie in 
different row and column). On the other hand, as we have seen in 
Proposition \ref{prop3.2}, $H=\{ {\rm lub}(\a,\b) : \a,\b\in 
\Gamma\}$. Since we can take $\phi_\a\in I$ except for finitely many 
$\a\in \Gamma$, we deduce that almost all elements in $H$ belong to 
$H(I)$. Thus $H(I)\cup\{\cero\}$ is also a semigroup of finite genus.
   \end{remark}
 
   \section{ The Main result}\label{s5}

Let $R$ be an $\bF$-algebra equipped with two well agreeing n-weights 
$\rho$ and $\sigma$.
   \begin{lemma}\label{lemma5.1} Let $f\in R^*.$ There exists $g\in 
\cM_\rho$ such that $fg\in\cM_\rho.$
   \end{lemma}
   \begin{proof} If $\rho(fg)=0$ for all $g\in \cM_\rho$ then, by 
property (N3), the same happens for all $g\in R^*$. Thus $I\subseteq 
\cU_{\rho}$, where $I=(f)$ is the ideal generated by $f$. This 
contradicts the fact that $H(I)\cup\{ \cero\}$ is a semigroup of 
finite genus.
   \end{proof} Define the map $\tilde{\rho}:R\rightarrow 
\bZ\cup\{-\infty\}$ as follows: $\tilde{\rho}(0):=-\infty$ and for 
$f\neq 0$,
  $$ \tilde{\rho}(f):=\min\{\rho(fg)-\rho(g): g\in\cM_\rho\}\, .
  $$ The following Lemma subsume some relevant properties of 
$\tilde{\rho}$.
   \begin{lemma}\label{lemma5.2}
  \begin{enumerate} \item $\tilde{\rho}(f)=\rho(fg)-\rho(g)$ for all 
$g\in\cM_{\rho}$ such that $fg\in\cM_\rho;$ \item If $f\in\cM_{\rho},$ 
then $\tilde{\rho}(f)=\rho(f)>0;$ if $f\in\cU_{\rho},$ then 
$\tilde{\rho}(f)\leq 0;$ \item $\tilde{\rho}(f)=0$ for all 
$f\in\bF^*;$ \item $\tilde{\rho}(fg)=\tilde{\rho}(f)+ 
\tilde{\rho}(g);$ \item $\tilde{\rho}(f+g)\leq\max\{\tilde{\rho}(f), 
\tilde{\rho}(g)\}.$
   \end{enumerate}
   \end{lemma}
   \begin{proof} (1) Let $g_1,g_2\in\cM_\rho$ such that 
$fg_1\in\cM_\rho$. Then 
$\rho(fg_1)+\rho(g_2)=\rho(fg_1g_2)\le\rho(fg_2)+\rho(g_1)$, and the 
right hand inequality is an equality when $fg_2\in\cM_\rho$. (2) and 
(3) are immediate. (4) By Lemma \ref{lemma5.1}, there exists 
$h\in\cM_\rho$ such that $fgh, gh\in\cM_\rho$. Then $\tilde{\rho}(fg)= 
(\rho(fgh)-\rho(gh))+(\rho(gh)-\rho(h))= \tilde{\rho}(f)+ 
\tilde{\rho}(g)$. (5) Let $h\in\cM_\rho$ such that $fh,gh\in\cM_\rho$. 
Then 
$\tilde{\rho}(f+g)\leq\rho((f+g)h)-\rho(h)\leq\max\{\tilde{\rho}(f), 
\tilde{\rho}(g)\}$.
   \end{proof} The additive inverse of $\tilde{\rho}$ can be extended 
to the whole field $\K$ in the usual way, so that we obtain a function 
$v_{\rho}$:
  \begin{equation*}\label{eq5.1} v_\rho(f/g):= \begin{cases}\infty & 
\text{if $f=0$,}\\ \tilde{\rho}(g)-\tilde{\rho}(f) & \text{if $f\neq 
0$.}
  \end{cases} \end{equation*} Properties in Lemma \ref{lemma5.2} 
implies the following.
   \begin{proposition}\label{prop5.1} The map $v_\rho$ is well defined 
and gives a discrete valuation of $\K$ over $\bF.$
   \end{proposition} Analogously we can define the valuation 
$v_\sigma$ associated to the n-weight $\sigma$. Denote by ${\mathbb 
P}(\K)$ the set of places of $\K$ over $\mathbb F$. For a place 
$S\in{\mathbb P}(\K)$, let $v_S$ and ${\mathcal O}_S$ be the 
corresponding valuation and valuation ring in $\K$. Set
  $$ 
{\mathcal S}(R):=\{S\in{\mathbb P}(\K): R\subseteq \cO_S\}\, .
  $$
   \begin{proposition}\label{prop5.2}{\rm (cf. 
\cite[p.2009]{matsumoto})} Let $P$ and $Q$ be the places of $\K$ 
corresponding to $v_\rho$ and $v_\sigma$ $($see Proposition 
\ref{prop5.1}$).$ Then
  $$ {\mathcal S}(R)=\bP(\K)\setminus\{P,Q\}\, .
  $$
   \end{proposition}
   \begin{proof} 
If $R\subseteq {\mathcal O}_P$ then $\cU_\rho=R$ hence 
$\cM_\rho=\emptyset$ and the semigroup $H(R)$ cannot have a finite 
genus. Thus $P,Q\not\in {\mathcal S}(R)$. Conversely, if ${\mathcal 
S}(R)\cup\{P,Q\}\neq \bP(\K)$, we can apply to ${\mathcal 
S}(R)\cup\{P,Q\}$ the Strong Approximation Theorem (see e.g. 
Stichtenoth \cite[I.6.4]{Stichtenoth}) to conclude that there exists a 
infinite sequence $(h_1,h_2,\ldots)$ of functions in $\K$ such that 
$v_\rho(h_i)=v_\sigma(h_i)=i$ and $v_S(h_i)\geq 0$ for each $S\in 
{\mathcal S}(R)$. In particular, $h_i\in \cap_{S\in {\mathcal 
S}(R)}\cO_S$ and this ring is precisely $\bar R$, the integral closure 
of $R$ in $\K$ (see e.g. \cite[III.2.6]{Stichtenoth}). The 
sequence $(h_1,h_2,\ldots)$ is $\bF$-linearly independent and 
contained in the $\bF$-vector space
   $$ 
W:=\{x\in \bar R: \text{$v_\rho(x)>0$ and $v_\sigma(x)>0$}\}\, .
   $$ 
As the n-orders $\rho$ and $\sigma$ are well agreeing, we have  
$W\cup R\subseteq \cU_\rho\cap\cU_\sigma=\bF$,  and thus
$W\cup R= \{0\}$. Then $\dim_\bF (W)\leq \dim_\bF(W+R)/R)\leq \dim_\bF 
(\bar R/R)$. But, according to
the Finiteness of Integral Closure Theorem 
(see e.g. \cite[Cor. 13.13]{eisenbud} or Zariski-Samuel 
\cite[Ch. V, Thm. 9]{zariski-samuel}), 
this last dimension is finite and we get a contradiction. 
   \end{proof}
Thus, we have proved the following.   
   \begin{theorem}\label{thm5.1} Let $R$ be an $\bF$-algebra admitting 
two well agreeing n-weights $\rho$ and $\sigma.$ Then
  \begin{enumerate} 
\item $R$ is an integral domain and its quotient 
field $\K$ is an algebraic function field of one variable over $\bF;$ 
\item There exist two places $P, Q\in {\mathbb P}(\K)$ such that 
$\rho$ and $\sigma$ are derived from the valuations associated to $P$ 
and $Q$ by the procedure stated in Example \ref{ex2.1}$;$ and 
\item $\bar R=\cap_{S\in \bP(\K)\setminus\{P,Q\}} \cO_S.$
  \end{enumerate}
  \end{theorem}
   \begin{remark}\label{rem5.1} Let $R$ be an integral domain 
$\bF$-algebra having Krull dimension 1, and let $\K$ be its field of 
quotients. Let $P,Q\in \bP(\K)$. By using the procedure of Example 
\ref{ex2.1}, the valuations at $P$ and $Q$ define two n-weights, 
$\rho$ and $\sigma$, over $R$.

Let us note that condition (3) in Theorem \ref{thm5.1} can be stated 
also as ${\mathcal S}(R)=\bP(\K)\setminus\{P,Q\}$. In this case 
$\cU_\rho\cap\cU_\sigma=\bF$. In fact, if $f\in R$ is such that 
$v_P(f)\geq 0$ and $v_Q(f)\geq 0$, then $f\in \cap_{S\in 
\bP(K)}\cO_S=\bF$ by \cite[III.2.6]{Stichtenoth}. Thus $\rho$ and 
$\sigma$ agree well iff $\# (\bN^2\setminus H(R))$ is finite. This 
observation leads to the following

{\em Question. Does condition (3) imply $\# (\bN^2\setminus 
H(R))<\infty$? }

So far we do not have an answer to this question.
  \end{remark}

   \end{document}